\magnification=1200
\noindent
LINEAR RESPONSE THEORY FOR DIFFEOMORPHISMS 
\medskip\noindent
WITH TANGENCIES OF STABLE AND UNSTABLE MANIFOLDS.
\medskip\noindent
[A CONTRIBUTION TO THE GALLAVOTTI-COHEN CHAOTIC HYPOTHESIS.]
\bigskip
\centerline{by David Ruelle\footnote{$\dagger$}{Math. Dept., Rutgers University, and 
IHES, 91440 Bures sur Yvette, France. email: ruelle@ihes.fr}.}
\bigskip\bigskip\bigskip\bigskip\noindent
	{\leftskip=1.8cm\rightskip=2cm{\sl Abstract:}
\medskip\noindent
This note presents a non-rigorous study of the linear response for an SRB (or `natural physical') measure $\rho$ of a diffeomorphism $f$ in the presence of tangencies of the stable and unstable manifolds of $\rho$.  We propose that generically, if $\rho$ has no zero Lyapunov exponent, if its stable dimension is sufficiently large (greater than 1/2 or perhaps 3/2) and if it is exponentially mixing in a suitable sense, then the following formal expression for the first derivative of $\rho(\phi)$ with respect to $f$ along $X$ is convergent:
$$	\Psi(z)=\sum_{n=0}^\infty z^n\int\rho(dx)\,X(x)\cdot\nabla_x(\phi\circ f^n)\qquad{\rm for}\qquad z=1   $$
This suggests that an SRB measure may exist for small perturbations of $f$, with weak differentiability.\par}
\vfill\eject
\bigskip\noindent
{\bf 1. Introduction.}
\medskip
	Let $f$ be a diffeomorphism of the compact manifold $M$ and $\ell$ a probability measure with smooth density.  Assume that $(f^n)^*$ applied to $\ell$ has a weak limit $\rho$ when $n\to+\infty$, so that $\rho$ is a `natural physical' or SRB measure.  If $f$ and $\rho$ depend on a real parameter $\alpha$, we say that we have {\it linear response} if $\alpha\to\rho_\alpha(\phi)$ is C$^1$ when $\phi:M\to{\bf R}$ is smooth.  This is a physically significant fact which is known to hold if $f_\alpha$ is uniformly hyperbolic (specifically if $\rho_\alpha$ is an SRB measure on a mixing Axiom A attractor for $f_\alpha$).
\medskip
	Uniform hyperbolicity is uncommon in physical `chaotic' systems (those for which $\rho$ has a Lyapunov exponent $>0$).  Usually chaos follows if $f$ effects some folding in the phase space $M$ and this folding results in tangencies of the stable and unstable manifolds of $\rho$, so that uniform hyperbolicity does not hold.  The conditional measures of $\rho$ on unstable manifolds have a density with a projection along the stable direction which (generically) has $1/${\it square root} singularity corresponding to tangencies.
\medskip
	If $f_\alpha$ and $\rho_\alpha$ depend on $\alpha$, linear response means thus that $\rho$ has a (weak) derivative with respect to $\alpha$.  This derivative can be formally computed, and the formal result may diverge because of the $1/$square root singularities mentioned above.  We shall however argue that if the unstable manifolds are piled up sufficiently densely, corresponding to a large stable dimension $d_S$ for $\rho_\alpha$ ($d_S>1/2$ or perhaps $>3/2$) this beats the $1/$square root singularities and linear response may hold.  We follow here the ideas of an earlier paper [15] with some corrections.
\medskip
	Note that $\rho$ is more likely to have a large stable dimension $d_S$ when dim$M$ increases (systems with many degrees of freedom).
\medskip
	The arguments presented here definitely do not constitute a rigorous proof but support a physically relevant conjecture in Section 6 below.  This conjecture agrees in particular with the fact (observed in computer experiments) that many nontrivial time evolutions behave as if they corresponded to uniformly hyperbolic dynamics (this is the {\it chaotic hypothesis} of Gallavotti-Cohen [2],[3]).  [A weaker and vaguer form of the conjecture is that a change of regime of linear response occurs when the unstable dimension of $\rho$ increases beyond 1/2].
\medskip
	We assume sufficient differentiability of $f$, say C$^3$.  Flows are not discussed in the present paper.
\bigskip\noindent
{\bf 2. Formal linear response formula.}
\medskip
	Let $M$ be a smooth compact manifold and $f:M\to M$ a diffeomorphism.  We denote by $\rho$ a `natural physical' measure on $M$.  We may take this to mean that $\rho$ is a weak limit
$$	\rho=\lim_{n\to\infty}(f^n)^*\ell\eqno{(1)}   $$
where $\ell$ is a probability measure on $M$ absolutely continuous with respect to {\it Lebesgue} (i.e., with respect to a Riemann volume).  In good cases (1) means that $\rho$ is an SRB measure.
\medskip
	SRB measures (Sinai, Ruelle, Bowen measures, see for instance [18]) were first defined for uniformly hyperbolic diffeomorphisms.  A modern mathematical definition of an SRB measure $\rho$ is that it is an $f$-invariant probability measure which satisfies either of the following equivalent conditions (Ledrappier-Strelcyn [5], Ledrappier-Young [6]):
\par\noindent
(i) the (Kolmogorov-Sinai) entropy of $(\rho,f)$ is the sum of the positive Lyapunov exponents,
\par\noindent
(ii) the conditional measures of $\rho$ on local unstable manifolds of $(\rho,f)$ are absolutely continuous with respect to Lebesgue on the unstable manifold.
\par\noindent
[We may take $\rho$ to be ergodic.  We shall define Lyapunov exponents and unstable manifolds in Section 4.]
\medskip
	If $X$ is a vector field on $M$, which we may think of as infinitesimal, we perturb $f:x\to fx$ to $\tilde f:x\to fx+X(fx)$.  We have to first order in $X$:
$$	\tilde f^kx=f^kx+\sum_{j=1}^k(T_{f^jx}f^{k-j})X(f^jx)   $$
where $T_xf$ denotes the tangent map to $f$ at the point $x$.
\medskip
	If the perturbation $\tilde f$ of $f$ replaces $\rho$ by an SRB measure $\rho+\delta\rho$ for $\tilde f$ we call $\delta\rho$ the linear response.  If $\phi$ is a smooth function on $M$ we have
$$	\phi(\tilde f^kx)=\phi(f^kx)+\phi'(f^kx)\sum_{j=1}^k(T_{f^jx}f^{k-j})X(f^jx)      $$
$$	=\phi(f^kx)+\sum_{j=1}^kX(f^jx)\cdot\nabla_{f^jx}(\phi\circ f^{k-j})      $$
hence
$$	\delta\rho(\phi)=\lim_{n\to\infty}{1\over n}
	\sum_{k=1}^n\sum_{j=1}^k\int\ell(dx)\,X(f^jx)\cdot\nabla_{f^jx}(\phi\circ f^{k-j})      $$
$$	=\lim_{n\to\infty}{1\over n}\sum_{k=1}^n\sum_{j=1}^k
	\int((f^j)^*\ell(dx))\,X(x)\cdot\nabla_x(\phi\circ f^{k-j})      $$
$$	=\lim_{n\to\infty}{1\over n}\sum_{i\ge0}\sum_{j=1}^{n-i}
	\int((f^j)^*\ell(dx))\,X(x)\cdot\nabla_x(\phi\circ f^i)      $$
If we interchange in the right-hand side $\lim_{n\to\infty}$ and $\sum_{i\ge0}$ (without a good mathematical justification), and use $\lim_{n\to\infty}{1\over n}\sum_{j=1}^{n-i}(f^j)^*\ell=\rho$ we obtain the formal linear response formula
$$	\delta\rho(\phi)=\sum_{n=0}^\infty\int\rho(dx)\,X(x)\cdot\nabla_x(\phi\circ f^n)=\Psi(1)\eqno{(2)}   $$
where we have defined the {\it susceptibility} function
$$	\Psi(z)=\sum_{n=0}^\infty z^n\int\rho(dx)\,X(x)\cdot\nabla_x(\phi\circ f^n)\eqno{(3)}   $$
\medskip
	Suppose that $f=f_\alpha\big|_{\alpha=0}$ where $\alpha$ is in a neighborhood of $0$ in ${\bf R}$ and $f_\alpha$ has SRB measure $\rho_\alpha$, we would like to know if
$$	[{d\over d\alpha}\int\rho_\alpha(dx)\,\phi(x)]\big|_{\alpha=0}=\Psi(1)\eqno{(4)}   $$
where $\Psi$ is defined with $X(fx)=df_\alpha x/d\alpha\big|_{\alpha=0}$.  In fact we shall be mostly concerned with the preliminary step of studying the convergence of $\Psi(1)$.
\bigskip\noindent
{\bf 3. Results on linear response.}
\medskip
	The main positive result is that if $\alpha\to f_\alpha$ is C$^3$ from a neighborhood of $0$ in ${\bf R}$ to C$^3$ diffeomorphisms of $M$, and $\rho_\alpha$ is SRB on a mixing Axiom A attractor for $f_\alpha$, then (4) holds if $\phi$ is C$^2$.  (For this result see Ruelle [16], and also Gou\"ezel-Liverani [4]).  We may here assume that $\rho_\alpha$ is defined by (1) with $f_\alpha$ replacing $f$, and $\ell$ has its support in a neighborhood of the compact attractor $A_\alpha$ such that $A_\alpha=\lim_{n\to\infty} f_\alpha^n\,{\rm supp}\,\ell$, where $f_\alpha$ satisfies uniform hyperbolicity assumptions on $A_\alpha$.
\medskip
	The above positive result is based on uniform hyperbolicity: the tangent bundle $T_AM$ over the attractor $A$ is the continuous direct sum of a contracting vector bundle $E^s$ and an expanding vector bundle $E^u$.  Write $X=X^s+X^u$ with $X^{s,u}\in E^{s,u}$, then (3) becomes
$$	\Psi(z)=\sum_{n=0}^\infty z^n\int\rho(dx)\,X^s(x)\cdot\nabla_x(\phi\circ f^n)
	-\sum_{n=0}^\infty z^n\int\rho(dx)\,({\rm div}^uX^u(x))(\phi(f^nx))\eqno{(5)}   $$
where div$^uX^u$ is the divergence of $X^u$ along a piece of unstable manifold containing $x$, with its natural measure.  The contracting term $\int\rho(dx)\,X^s(x)\cdot\nabla_x(\phi\circ f^n)$ tends to 0 exponentially with respect to $n$.  The expanding term $\int\rho(dx)\,({\rm div}^uX^u(x))(\phi(f^nx))$ is a correlation function of ${\rm div}^uX^u(x)$ and $\phi$ which also tends to 0 exponentially with respect to $n$.  Therefore the radius of convergence of $\Psi$ is $>1$ and $\Psi(1)$ is well defined.
\medskip
	Another result holds if $f$ preserves a smooth volume $\rho$ on $M$ and if $(f,\rho)$ is exponentially mixing.  Then (3) can be rewritten as
$$	\Psi(z)=-\sum_{n=0}^\infty\int\rho(dx)({\rm div}X(x))(\phi(f^nx))   $$
and the radius of convergence of $\Psi$ is $>1$.
\medskip
	In computer studies, nontrivial dynamics of moderate dimension can sometimes lead to very long characteristic times, or sometimes to systems that behave like uniformly hyperbolic systems (chaotic hypothesis [2], [3], which remains to be explained).
\medskip
	There are a number of other results, positive and negative, for which we refer to Baladi [1].  In particular for interval maps it is known that $\alpha\to\rho_\alpha$ is often discontinuous [17]; in fact $X^u(x)=X(x)$ often has singularities $\sim1/\sqrt{\pm(x-x_j)}$ so that ${\rm div}^uX(x)$ does not make sense as a function.\footnote{*}{We use $\sim$ to denote proportionality, while $\approx$ will denote approximate equality.}
\medskip
	There is an important literature on stable-unstable tangencies typified by Newhouse's infinitely many sinks (see [8], [9], [10], [11]).  This literature refers to dynamically unstable situations.  Here we try to investigate the opposite stable situation where SRB measures change differentiably under a small change of the diffeomorphism $f$.
\bigskip\noindent
{\bf 4. Stable and unstable manifolds.}
\medskip
	Let $\rho$ be an ergodic measure on an attractor $A$ for the diffeomorphism $f$ of the manifold $M$.  In this Section $\rho$ need not be SRB but we assume that $\rho$ is hyperbolic namely that all Lyapunov exponents (defined below) of $\rho$ are different from 0.  We describe now briefly the theory of stable and unstable manifolds (Pesin [12], [13] in a version due to Ruelle [14]).
\medskip
	Define $T_x^n=T_xf^n$, where $T_xf:T_xM\to T_{fx}M$ is the tangent map to $f$ at $x\in M$.  For $\rho$-almost all $x$ one can show that the following limit exists:
$$	\lim_{n\to\infty}(T_x^{n*}T_x^n)^{1/2n}=\Lambda_x\eqno{(6)}   $$
and its eigenvalues $\exp\lambda^{(1)}<\dots<\exp\lambda^{(\kappa)}$ with multiplicities $m^{(1)},\dots,m^{(\kappa)}$ are $x$-independent for $\rho$-almost all $x$ (note that we have $\kappa\le$ dim$M$).  The $\lambda^{(r)}$ are called the Lyapunov exponents of $\rho$.  Let $U_x^{(1)},\dots,U_x^{(\kappa)}$ be the eigenspaces of $\Lambda_x$.  Writing $V_x^{(0)}=\{0\}$ and $V_x^{(r)}=U_x^{(1)}\oplus\dots\oplus U_x^{(r)}$, we have
$$	\lim_{n\to\infty}{1\over n}\log||T_x^nu||=\lambda_x^{(r)}\qquad
	{\rm when}\qquad u\in V_x^{(r)}\backslash V_x^{(r-1)}\eqno{(7)}   $$
If $\lambda^{(s)}$  is the largest negative Lyapunov exponent the space $V_x^{(s)}=E_x^s$ is the stable subspace at $x$ and $(E_x^s)\big|_{x\in A}=E^s$ the contracting vector bundle.  Replacing $f$ by $f^{-1}$ replaces $E_x^s$ by the unstable subspace $E_x^u$ and $(E_x^u)\big|_{x\in A}=E^u$ is the expanding vector bundle.
\medskip
	If $0<a<-\lambda^{(s)}$ and $b>0$ we define a local exponential stable manifold
$$	{\cal V}_x^s=\{y\in M:d(f^ny,f^nx)<b.e^{-na}\quad\hbox{for all}\quad n\ge0\}\eqno{(8)}   $$
This stable manifold is a nonlinear version of the stable subspace $E_x^s$ as follows from (7).  The definition (8) implies that $T_{f^{(n-1)}y}f$ is exponentially close to $T_{f^{(n-1)}x}f$ for $n\to\infty$.  From this one obtains a perturbation theorem: if $b$ is sufficiently small and $y\in{\cal V}_x^s$, then $\Lambda_y$ has the same eigenvalues as $\Lambda_x$ including multiplicity and there is $B>0$ such that $||V_y^{(r)}-V_x^{(r)}||\le Bd(y,x)$.  In fact $E_y^s$ is the tangent space to ${\cal V}_y^s$ at $y\in{\cal V}_x$ and depends continuously on $y$.  Therefore ${\cal V}_x^s$ is a differentiable manifold.  One can deal with higher derivatives as with $T_xf$ so that the stable manifolds are as smooth as $f$.
\medskip
	Similarly for the unstable manifolds ${\cal V}_x^u$.  Note furthermore that these unstable manifolds are contained in the attractor $A$.
\bigskip\noindent
{\bf 5. Stable-unstable tangencies.}
\medskip
	A simple feature of the diffeomorphism $f$ which can lead to positive Lyapunov exponents (chaotic behavior) and exponential decay of correlations is when $f$ folds the attractor $A$.  This folding causes tangencies of the stable and unstable manifolds, which we now want to study.  The set of points of tangency of a stable and an unstable manifold has $\rho$-measure 0 in $M$ since $\rho$ has no zero Lyapunov exponent.  One can nevertheless deal with stable-unstable tangencies as follows.
\medskip
	In the framework of Section 4 the stable manifolds ${\cal V}_x^s$ and the unstable manifolds ${\cal V}_x^u$, defined $\rho$-almost everywhere, depend measurably on $x$ together with their tangential derivatives.  Let the open neighborhoods ${\cal N}$ of a point of $M$ have measure $\rho({\cal N})>0$ and write again ${\cal V}_x^s$ and ${\cal V}_x^u$ for the restrictions of these manifolds to ${\cal N}$.  Using Lusin's theorem we can choose ${\cal N}$ and a family $({\cal V}_x^s)_{x\in S}$ with compact $S\subset{\cal N}$, where each ${\cal V}_x^s$ contains only one $x\in S$, such that the ${\cal V}_x^s$ and tangential derivatives depend continuously on $x$ and $\rho({\cal N}\backslash\cup_{x\in S}{\cal V}_x^s)<\epsilon$ for small $\epsilon$.  We can also choose a similar family $({\cal V}_{y}^u)_{y\in U}$ of unstable manifolds.
\medskip
	In brief we can, locally and up to a set of small $\rho$-measure, consider that the stable manifolds form a family (a pile) $({\cal V}_x^s)_{x\in S}$ continuous together with their tangential derivatives.  A pile $({\cal V}_{y}^u)_{y\in U}$ is similarly defined for unstable manifolds.  Since we have control of the second derivatives of the stable and unstable manifolds we can impose that the intersections ${\cal V}_x^s\cap{\cal V}_{y}^u$ are interior to ${\cal V}_x^s$ and ${\cal V}_{y}^u$ and each is either empty, or consists of two points, or of one regular tangency point of ${\cal V}_x^s$, ${\cal V}_{y}^u$.
\medskip
	We shall take coordinates such that the stable manifolds $({\cal V}_x^s)_{x\in S}$ are roughly parallel and the unstable manifolds $({\cal V}_{y}^u)_{y\in U}$ are folded so that the intersections ${\cal V}_x^s\cap{\cal V}_{y}^u$ are as above.  Note that the folding of unstable manifolds which gives rise to tangencies gives folds dense in the support of $\rho$, and this explains why part of ${\cal N}$ has to be excluded from the pile $({\cal V}_{y}^u)_{y\in U}$.  Similarly for stable manifolds.
\medskip
	We choose now a manifold ${\cal W}$ with dim${\cal W}=$ dim${\cal V}^u$ such that ${\cal W}$ is `parallel to the unstable fold' namely transversal to the ${\cal V}_x^s$, and we project ${\cal N}$ on ${\cal W}$ along the ${\cal V}_x^s$.  This projection $\varpi$ is `almost' absolutely continuous from each ${\cal V}_y^u$ with its natural measure to ${\cal W}$ with Lebesgue measure (some Riemann measure).  The corresponding Jacobian $J(x)$ is `almost' H\"older continuous with H\"older exponent close to 1 except along the fold (the stable-unstable tangency points of ${\cal V}_y^v$).  The Jacobian (generically) has a $1/$square root singularity transversal to the projection of the fold on ${\cal W}$.
\bigskip\noindent
{\bf 6. A conjecture on linear response in the presence of stable-unstable tangencies.}
\medskip
	In the presence of regular tangencies of stable-unstable manifolds (a condition on second-order derivatives for the piles $({\cal V}_x^s)_{x\in S}$ and $({\cal V}_{y}^u)_{y\in U}$), and with a genericity assumption on the projection along the stable direction, we expect the following:
\medskip
	 {\sl If the SRB measure $\rho$ for the diffeomorphism $f$ has no zero Lyapunov exponent, if its stable dimension is sufficiently large (greater than 1/2 or perhaps 3/2) and if it is exponentially mixing in a suitable sense, then the following formal expression for the first derivative of SRB with respect to $f$ along $X$ evaluated at $\phi$ is convergent:
$$	\Psi(z)=\sum_{n=0}^\infty z^n\int\rho(dx)\,X(x)\cdot\nabla_x(\phi\circ f^n)\qquad{\rm for}\qquad z=1.   $$ }
\medskip
	The nonrigorous nature of the following study of $\Psi$ results in part from the fact that we shall (in the spirit of [15]) ignore the small measure sets ${\cal N}\backslash\cup_{x\in S}{\cal V}_x^s$ and ${\cal N}\backslash\cup_{y\in U}{\cal V}_y^u$.
\bigskip\noindent
{\bf 7. Estimating the radius of convergence of $\Psi$.}
\medskip
	We consider the contribution $\Psi_{\cal N}(z)$ of $\rho\big|{\cal N}$ to the susceptibility:
$$	\Psi_{\cal N}(z)=\sum_{n=0}^\infty z^n\int_{\cal N}\rho(d\xi)\,X(\xi)\cdot\nabla_\xi(\phi\circ f^n)\eqno{(9)}   $$
in order to estimate how the fold in ${\cal N}$ influences the radius of convergence of $\Psi$.
\medskip
Write $\rho(d\xi)=\int\sigma(dx)\rho_x(d\xi)$ where $\rho_x$ is the natural measure on the unstable manifold ${\cal V}_x^u$ and $\sigma(dx)$ is a `stable' measure carried by a set of Hausdorff dimension d$_S$.  The projection $\varpi[X(\xi)\rho_x(d\xi)]$ has a `almost' H\"older density on ${\cal W}$ (H\"older exponent close to 1 if ${\cal N}$ is small) with a $1/$square root singularity transversal to the projection of the fold of ${\cal V}_x^u$ on ${\cal W}$.  [We have used here the philosophy that we can ignore the small-measure sets ${\cal N}\backslash\cup_{x\in S}{\cal V}_x^s$ and ${\cal N}\backslash\cup_{y\in U}{\cal V}_y^u$].
\medskip
	We parametrize ${\cal W}$ by variables $w,\theta$ with $\theta\in{\bf R}$ such that the lines $\ell_w$ given by $w=$ constant (and parametrized by $\theta$) are transversal to the folds.
\medskip
	The projection $\varpi[X(\xi)\rho(d\xi)]$ restricted to $\ell_w$ contains a factor due to $X(\xi)$ and we can estimate its density with respect to $\theta$ as
$$	\Delta(\theta)=\int\sigma(dx)\Delta_x(\theta)\eqno{(10)}   $$
Here $\Delta_x(\theta)$ is the projection $\varpi[X(\xi)\rho_x(d\xi)]$ restricted to $\ell_w$ so that $\Delta_x(\theta)\sim1/\sqrt{\pm(\theta-\theta_x)}$ where $\theta_x$ depends on $w$ [the presence of the square root rather than another power is a genericity assumption as mentioned in Section 6].
\medskip
	Define $\psi(\tau)=\int_{\{\theta_x<\tau\}}\sigma(dx)$ or $\psi(\tau)=-\int_{\{\theta_x>\tau\}}\sigma(dx)$.  The measure $d(\psi(\tau))$ is a 1-dimensional projection of the $d_S$-dimensional measure $\sigma(dx)$.  Using Marstrand's theorem [7] we assume that $d(\psi(\tau))$ has dimension $\bar d_S$=min$\{1,d_S\}$ [this is again a genericity assumption as mentioned in Section 6].  Therefore $|\psi(\theta)-\psi(\tau)|<|\theta-\tau|^{\bar d_S}$ up to a multiplicative constant.  We can estimate $\Delta(\theta)$ in terms of
$$	\int\sigma(dx)\,\Delta_x(\theta)\sim\int_{\tau<\theta}{d(\psi(\tau))\over\sqrt{\theta-\tau}}
	+\int_{\tau>\theta}{d(\psi(\tau))\over\sqrt{\tau-\theta}}
	<2\bar d_S\int_0^Tt^{\bar d_S-3/2}\,dt   $$
with $T$ depending on the support of $\rho$.
\medskip
	If $d_S>1/2$  we show now that $\Delta$ is H\"older continuous [the proof of this fact in [15] is incorrect].  For $\delta>0$ we can estimate $\Delta(\theta+\delta)-\Delta(\theta)$ in terms of
$$ 	\int_{-T}^{\theta+\delta}{d(\psi(\tau))\over\sqrt{\theta+\delta-\tau}}
	-\int_{-T}^\theta{d(\psi(\tau))\over\sqrt{\theta-\tau}}\qquad{\rm and}\qquad
	\int_\theta^T{d(\psi(\tau))\over\sqrt{\tau-\theta}}
	-\int_{\theta+\delta}^T{d(\psi(\tau))\over\sqrt{\tau-\delta-\theta}}\eqno{(11)}   $$
The first term of (11) is
$$	=\int_\theta^{\theta+\delta}{d(\psi(\tau))\over\sqrt{\theta+\delta-\tau}}
	+\int_{\theta-\delta}^\theta{d(\psi(\tau))\over\sqrt{\theta+\delta-\tau}}
	-\int_{\theta-\delta}^\theta{d(\psi(\tau))\over\sqrt{\theta-\tau}}\eqno{(12)}   $$
$$    +\int_{-T}^{\theta-\delta}\Big({1\over\sqrt{\theta+\delta-\tau}}-{1\over\sqrt{\theta-\tau}}\Big)d(\psi(\tau))   $$
With $t=\theta+\delta-\tau$ we have
$$	\Big|\int_\theta^{\theta+\delta}{d(\psi(\tau))\over\sqrt{\theta+\delta-\tau}}\Big|
	\sim\int_0^\delta{|d(\psi(t))|\over\sqrt{t}}<\bar d_S\int_0^\delta{t^{\bar d_S}dt\over t^{3/2}}
	=\Big({\bar d_S\over \bar d_S-1/2}\Big)\delta^{\bar d_S-1/2}\eqno{(13)}   $$
and similar estimates for the next two terms of (12).  Furthermore
$$	\Big|\int_{-T}^{\theta-\delta}\Big({1\over\sqrt{\theta+\delta-\tau}}-{1\over\sqrt{\theta-\tau}}\Big)d(\psi(\tau)\Big|
	<{\delta\over2}\int_{-T}^{\theta-\delta}{|d(\psi(\tau))|(\theta-\tau)^{-3/2}}   $$
$$	<\bar d_S{\delta\over2}\int_{-T}^{\theta-\delta}(\theta-\delta-\tau)^{\bar d_S-1}(\theta-\tau)^{-3/2}\,d\tau
\approx\bar d_S{3\delta\over4}\int_{-T}^{\theta-\delta}(\theta-\delta-\tau)^{\bar d_S-1}(\theta-\tau)^{-5/2}\,d\tau   $$
$$	<\bar d_S{3\delta\over4}\int_\delta^{\theta+T}(t-\delta)^{\bar d_S}t^{-5/2}\,dt
	<\bar d_S{3\delta\over4}\int_\delta^{\theta+T}t^{\bar d_S-5/2}\,dt
	\sim\delta[-t^{\bar d_D-3/2}]\big|_\delta^{\theta+T}\approx\delta^{\bar d_D-1/2}\eqno{(14)}   $$
From (12--14) and similar results for the second term of (11) we see that $\Delta$ is $(\bar d_D-1/2)$-H\"older continuous.
\medskip
	We note that the H\"older continuity of $\Delta$ implies that $\int_{\cal N}\rho(d\xi)\,X(\xi)\cdot\nabla_\xi(\phi\circ f^n)\to0$ exponentially when $n\to\infty$ if the correlations for $\rho$ decay exponentially in a suitable sense.  In this case $\Psi_{\cal N}(1)$ is thus well defined.
\medskip
	Writing $d_S$ as a sum of partial dimensions one may conjecture that $\psi$ is $d_S$-H\"older, and therefore differentiable if $d_S>1$.  The derivative $\Delta'$ of $\Delta$ satisfies then
$$	\Delta'(\theta)\sim\int_{-T}^\theta{d\psi'(\tau)\over\sqrt{\theta-\tau}}
	+\int_\theta^T{d\psi'(\tau)\over\sqrt{\tau-\theta}}   $$
The same argument as above shows that $\Delta'$ is $({\rm min}(2,d_S)-3/2)$-H\"older if $d_S>3/2$.  The divergence ${\cal D}(\eta)$ of the projection $\varpi[X(\xi)\rho(d\xi)]$ on ${\cal W}$ is thus also $({\rm min}(2,d_S)-3/2)$-H\"older continuous (this results from the arbitrariness of the choice of the line $\ell_w$ if dim${\cal W}>1$).  Returning to (9) we see that
$$	\int_{\cal N}\rho(d\xi)\,X(\xi)\cdot\nabla_\xi(\phi\circ f^n)\sim-\int d\eta\,{\cal D}(\eta)\phi(f^n\eta)    $$
is a correlation function of H\"older functions which we assume to decay exponentially.  Therefore the radius of convergence of $\Psi_{\cal N}(z)$ is $>1$.  Since $\Psi_{\cal N}(z)$ is the most singular part of $\Psi(z)$, we find that the radius of convergence of $\Psi(z)$ is $>1$ and that the formal expression of the derivative (4) converges.
\bigskip\noindent
{\bf 8. Conclusions.}
\medskip
	To summarize we have shown that if the stable dimension $d_S$ of $\rho$ is $>1/2$ then $\varpi[X(\xi)\rho(d\xi)]$ is `approximately' H\"older, and if $d_S>3/2$ then $\varpi[(\nabla_\xi\cdot X)\rho(d\xi)]$ is `approximately' H\"older.  These H\"older conditions and a suitable property of exponential decay of correlations for $\rho$ imply that $\Psi_{\cal N}(1)$ and $\Psi(1)$ are well defined.  Therefore the derivative w.r.t. $f$ of the SRB measure $\rho$ makes sense generically in the presence of stable-unstable tangencies if $d_S$ is sufficiently large and there is a suitable exponential decay of correlations.
\bigskip\noindent
{\bf References.}
\medskip
[1] V. Baladi  ``Linear response, or else.''  arXiv:1408.2937(2014).
\medskip
[2] G. Gallavotti and E.G.D. Cohen  ``Dynamical ensembles in nonequilibrium statistical mechanics.''  Phys. Rev. Lett. {\bf 74},2694-2697(1995).
\medskip
[3] G. Gallavotti and E.G.D. Cohen  ``Dynamical ensembles in stationary states.''  J. Stat. Phys. {\bf 80},931-970(1995).
\medskip
[4] S. Gou\" ezel and C. Liverani  ``Banach spaces adapted to Anosov systems.''  Ergodic Theory
Dynam. Systems {\bf 26},189-217(2006).
\medskip
[5] F. Ledrappier and J.M. Strelcyn  ``A proof of the estimation from below in Pesin's
entropy formula.''  Ergodic Theory Dynam. Systems {\bf 2},203-219(1982).
\medskip
[6] F.Ledrappier and L.S.Young  ``The metric entropy of diffeomorphisms: I. Characterization of measures satisfying Pesin's formula, II. Relations between entropy, exponents and dimension.''  Ann. of Math. {\bf 122},509-539,540-574(1985).
\medskip
[7] J.M. Marstrand  ``Some fundamental geometrical properties of plane sets of fractional dimension.''  Proc. London Math. Soc.(3) {\bf 4}257-302(1954).
\medskip
[8] S.E. Newhouse  ``Diffeomorphisms with infinitely many sinks.''  Topology {\bf 13(1)},9-18(1974).
\medskip
[9] S.E. Newhouse  ``The Abundance of wild hyperbolic sets and nonsmooth stable sets for diffeomorphisms.''  IHES Publications Math\'ematiques {\bf 50},101-151(1979).
\medskip
[10] J. Palis and F. Takens.  {\it Hyperbolicity and Sensitive-Chaotic Dynamics at Homoclinic Bifurcations, Fractal Dimensions and Infinitely Many Attractors.}  Cambridge Univ. Press, 1993; Second Edition, 1994.
\medskip
[11] J. Palis and J.-Chr. Yoccoz  ``Non-uniformly hyperbolic horseshoes arising from bifurcations of heteroclinic cycles.''  IHES Publications Math\'ematiques {\bf 110},1-217(2010).
\medskip
[12] Ya.B.Pesin.  ``Invariant manifold families which correspond to non-vanishing characteristic exponents.''  Izv. Akad. Nauk SSSR Ser. Mat. {\bf 40},No 6,1332-1379(1976).  English translation: Math. USSR Izv. {\bf 10},No 6,1261-1305(1976).
\medskip
[13] Ya.B.Pesin.  ``Lyapunov characteristic exponents and smooth ergodic theory.''  Uspehi Mat. Nauk {\bf 32},No 4,55-112(1977).  English translation:  Russian Math. Surveys. {\bf 32},No 4,55-114(1977).
\medskip
[14] D.Ruelle.  ``Ergodic theory of differentiable dynamical systems.''  IHES Publications Math\'ematiques {\bf 50},27-58(1979).
\medskip
[15] D. Ruelle  ``Singularities of the susceptibility of an SRB measure in the presence of stable-unstable tangencies.''  Phil. Trans. R. Soc. A {\bf 369},482-493(2011).
\medskip
[16] D. Ruelle  ``Differentiation of SRB states.''  Comm. Math. Phys. {\bf 187},227-241(1997); {\bf 234},185-190(2003).
\medskip
[17] M. Yakobson  ``Absolutely continuous invariant measures for one-parameter families of one-dimensional maps.''  Commun. Math. Phys. {\bf 81},39-88(1981).
\medskip
[18] L.-S. Young  ``What are SRB measures, and which dynamical systems have them?''  J. Statist. Phys. {\bf 108},733-754(2002).
\end